\documentclass[letterpaper,12pt]{article}

\usepackage{amssymb}
\usepackage{amsmath, amscd}
\usepackage{graphicx}

\def\ni{\noindent}

\def\nn{\mathbb{N}}
\def\qq{\mathbb{Q}}
\def\rr{\mathbb{R}}
\def\ff{\mathbb{F}}
\def\cc{\mathbb{C}}
\def\zz{\mathbb{Z}}
\def\gg{\mathbb{G}}
\def\pp{\mathbb{P}}

\begin{document}

\title{Zeta Functions of Projective Toric Hypersurfaces over Finite Fields}
\author{Chiu Fai Wong}
\date{ }

\maketitle

\begin{abstract}
I give a formula for the zeta function of a projective toric hypersurface over a finite field and estimate its Newton polygon. As an application this formula allows us to compute the exact number of rational points on the families of Calabi-Yau manifolds in Mirror Symmetry. 
\end{abstract}
\ni

\vspace{2mm}
\ni

\section{Introduction}
\

Let $\ff_q$ be a finite field of $q$ elements where $q$ is a power of a prime $p$. Let $f= \sum_{j \in J}a_jx^j \in \ff_q[x_1^{\pm 1}, \ldots , x_n^{\pm 1}]$ be a Laurent polynomial over $\ff_q$, where the sum is over a finite subset $J$ of $\zz^n$. We assume $a_j \neq 0$ for all $j \in J$. Let $U_f$ be the affine hypersurface in the $n$-torus $\gg_m^n:=(\ff_q^*)^n$ defined by $\{ f=0 \}$. We can construct a projective toric variety $\pp_{\bigtriangleup}$ from the polytope $\bigtriangleup (f)$, the convex hull in $\rr^n$ of the lattice points which occur as exponents of the monomials of $f$ together with the origin. Assume $\dim \bigtriangleup (f) =n$. According to the theory of toric varieties, $\pp_{\bigtriangleup}$ is a compactification of the $n$-dimensional torus $\gg_m^{n}$ by algebraic tori $\gg_m^{\dim \sigma}$ of smaller dimensions. Each $\gg_m^{\dim \sigma}$ corresponds to a nonempty face $\sigma \subseteq \bigtriangleup (f)$ of the polytope $\bigtriangleup (f)$ \cite{Dan}.\\ 

This paper gives a formula for the zeta function of the projective closure $\overline{U_f}$ of the affine toric hypersurface $U_f$ in $\pp_{\bigtriangleup}$ over a finite field $\ff_q$ when $\bigtriangleup (f)$ is a simplex. This result extends Dwork's work on a smooth, projective hypersurface in $\pp^n$.\\

It is natural to begin with the zeta function on an affine hypersurface $U_f$ on $\gg_m^n$ defined by $f$, and extend via the usual toric decomposition of $\pp_{\bigtriangleup}$ to the zeta function of the toric projective closure $\overline{U_f}$. The zeta function of the toric projective closure of an affine toric hypersurface is determined by the product of the $L$-functions associated to each face $\sigma$ of $\bigtriangleup (f)$.\\ 

Our methods are $p$-adic and are based on the work of Dwork \cite{Dwork}, \cite{Dwork2}. The main theorem (Theorem 4.7) is the extension of Dwork's cohomology theory from smooth, projective hypersurfaces in characteristic $p$ to $\bigtriangleup$-regular, toric projective hypersurfaces with simplex $\bigtriangleup (f)$. Given a Laurent polynomial $f$ on $\gg_m^n$, we construct a complex of a $p$-adic Banach space $B$ on which a Frobenius operator $\alpha$ acts. In fact, the Banach space $B$ is generated by the monomials of $x_0f$. The exponents of the monomials of $B$ can be identified with the lattice points in the monoid $M(x_0f)$ generated by the vertices of $1 \times \bigtriangleup (f)$ and the origin in $\rr^{n+1}$. Dwork trace formula implies that $L(x_0f, t)^{(-1)^n}$ is the alternating product of the characteristic polynomial of the Frobenius operator $\alpha$ which acts on the $p$-adic Koszul complex induced by $\bigtriangleup (f)$. We have the same result for each face $\sigma$ of $\bigtriangleup (f)$.\\

The main observation is that when $\bigtriangleup (f)$ is a simplex, there is a bijection between nonempty faces of $\bigtriangleup (f)$ and nonempty subsets of the vertices of $\bigtriangleup (f)$. The product of the $L$-functions
is the alternating product of the characteristic polynomials of Frobenius operators corresponding to each nonempty face. An inclusion-exclusion argument shows that the alternating product is determined by the subcomplex corresponding to the interior cone of the monoid $M(x_0 f)$. Passing to homology and using the acyclicity of $\bigtriangleup$-regularity, we have
\begin{eqnarray}
Z(\overline{U_f}, qt) = \frac{P(t)^{(-1)^n}}{(1-qt) \cdots (1-q^nt)}
\end{eqnarray}
where $P(t)$ is a polynomial of degree $(-1)^{n+1} \sum_{\sigma \subseteq \bigtriangleup (f)} (-1)^{\dim \sigma} ( \dim \sigma )! \mathrm{Vol} (\sigma)$.\\ 

As a consequence, we are able to bound below the Newton polygon of $P(t)$ by its Hodge polygon and their endpoints coincide. If $\bigtriangleup (f)$ is not a simplex, we do not expect that equation (1) holds in general.\\ 

Equation (1) can be used to compute the number of rational points on a mirror pair of Calabi-Yau hypersufaces in arithmetic mirror symmetry when $\bigtriangleup (f)$ is a reflexive simplex.\\

This paper is organized as follows. In section 2 and 3, we recall some basic technique developed from Dwork \cite{Dwork}, \cite{Dwork2}, Adolphson and Sperber \cite{AS1}, \cite{AS2}. In section 4, we prove an analogous formula of the zeta function of projective toric hypersurfaces for simplex case. In section 5, we give a sharp lower bound for Newton polygon of the zeta function. In section 6, we compute the zeta functions of the families of Calabi-Yau hypersurfaces.\\ 

I would like to express my hearty thanks to my advisor, Professor Daqing Wan, for his patient guidance, continual encouragement and help throughout the period of my doctoral studies and in the preparation of this paper.\\

\section{$p$-adic theory}
\

Let $p$ be a prime number, $\qq_p$ the field of $p$-adic numbers and $\zz_p$ the ring of $p$-adic integers. Let $\Omega$ be the completion of an algebraic closure of $\qq_p$. Let $\ff_q$ be the finite field of $q=p^a$ elements, $\ff_{q^k}$ the extension of $\ff_q$ of degree $k$, and denote by $K$ the unramified extension of $\qq_p$ in $\Omega$ of degree $a$. Let $\pi \in \Omega$ satisfy $\pi^{p-1}=-p$. Then $\Omega_1= \qq_p(\pi)$ is a totally ramified extension of $\qq_p$ of degree $p-1$, in fact, $\Omega_1= \qq_p(\zeta_p)$, where $\zeta_p$ is a primitive $p$-th root of unity. Let $\Omega_0=K(\pi)$ be the compositum of $\Omega_1$ and $K$. Then $\Omega_0$ is an unramified extension of $\Omega_1$ of degree $a$. The residue class fields of $\Omega_0$ and $K$ are both $\ff_q$, and the residue fields of $\Omega_1$ and $\qq_p$ are both $\ff_p$. The Frobenius automorphism $x \mapsto x^p$ of $Gal(\ff_q / \ff_p)$ lifts to a generator $\tau$ of $Gal(\Omega_0 / \Omega_1) ( \cong Gal(K/ \qq_p))$ which is extended to $\Omega_0$ by requiring $\tau (\pi)= \pi$. If $\zeta$ is a $(q-1)$-th root of unity in $\Omega_0$, then $\tau (\zeta)= \zeta^p$. Denote by ``ord" the additive valuation on $\Omega$ normalized by ord$p=1$, and denote by ``ord$_q$" the additive valuation normalized by ord$_qq=1$.\\         

Let $f$ be a Laurent polynomial on the torus $\gg_m^n:=(\ff_q^*)^n$ and write 
$$f= \sum_{j \in J}a_jx^j \in \ff_q[x_1^{\pm 1}, \ldots , x_n^{\pm 1}],$$
where the sum is over a finite subset $J$ of $\zz^n$. We may assume $a_j \neq 0$ for all $j \in J$. Each Laurent polynomial $f$ defines an affine toric hypersurface 
$$U_f= \{ (x_1, \ldots , x_n) \in \gg_m^n| f(x_1, \ldots , x_n)=0  \}$$
in $\gg_m^n$.\\

Define the polytope $\bigtriangleup:= \bigtriangleup (f)$ of $f$ to be the convex hull in $\rr^n$ of the lattice points which occur as exponents of the monomials of $f$ together with the origin. We can construct a toric variety $\pp_{\bigtriangleup}$ from the polytope $\bigtriangleup (f)$ as follows:\\

\noindent $\bf{Definition.}$ The polytope ring $S_{\bigtriangleup}$ of a polytope $\bigtriangleup$ is defined by 
$$S_{\bigtriangleup}= \ff_q[x_0^rx^m:m=(m_1, \ldots ,m_n) \in r \bigtriangleup \cap \zz^n, r \geq 0],$$
the subalgebra of $\ff_q[x_0, x_1^{\pm 1}, \ldots , x_n^{\pm 1}]$ generated as a $\ff_q$-vector space by 1 and all monomials $x_0^rx^m=x_0^rx_1^{m_1} \cdots x_n^{m_n}$ where $r \in \nn$ such that the rational point $m/r=(m_1/r, \ldots ,m_n/r)$ belongs to $\bigtriangleup$.\\

The standard grading of the polynomial ring $\ff_q[x_0, x_1^{\pm 1}, \ldots , x_n^{\pm 1}]$ in $x_0$ induces the grading of $S_{\bigtriangleup}$:
$$\deg (x_0^rx^m)=r$$
and the graded commutative $\ff_q$-algebra
$$S_{\bigtriangleup} = \bigoplus_{r \geq 0} S_{\bigtriangleup}^r,$$
where $S_{\bigtriangleup}^r$ is a finite dimensional $\ff_q$-vector space with the basis $x_0^rx^m$ $(m \in r \bigtriangleup \cap \zz^n)$.\\ 

Let $\bigtriangleup (x_0f) \subseteq \rr^{n+1}$ be the convex hull of the origin in $\rr^{n+1}$ and $(1, v)$ where $v \in \bigtriangleup (f)$. Define the monoid $M(x_0f) \subseteq \zz^{n+1}$ to be the lattice points generated by $\bigtriangleup (x_0f)$. It can be identified with the exponents of the monomials of the polytope ring $S_{\bigtriangleup}$.\\

We now can associate with the polytope $\bigtriangleup \subseteq \rr^n$ the projective algebraic variety $\pp_{\bigtriangleup}$.\\

\noindent $\bf{Definition.}$ The algebraic variety 
$$\pp_{\bigtriangleup}= \mathrm{Proj}(S_{\bigtriangleup})$$
is called the toric variety associated with polytope $\bigtriangleup$.\\

Given a face $\sigma$ of $\bigtriangleup (f)$, we define $\dim \sigma$ to be the dimension of the smallest subspace of $\rr^n$ containing $\sigma$. Assume $\dim \bigtriangleup (f) =n$. Define $\mathrm{Vol} (f)$ to be the volume of $\bigtriangleup (f)$ with respect to Lebesgue measure on $\rr^n$.\\

Toric vaiety $\pp_{\bigtriangleup}$ is a compactification of the $n$-dimensional torus $\gg_m^n$ by algebraic tori $\gg_m^{\dim \sigma}$ of smaller dimensions. Each $\gg_m^{\dim \sigma}$ corresponds to a nonempty face $\sigma \subseteq \bigtriangleup$ of the polytope $\bigtriangleup$ \cite{Dan}. That is, we have the so called toric decomposition 
$$\pp_{\bigtriangleup}= \bigsqcup_{\emptyset \neq \sigma \subseteq \bigtriangleup} \gg_m^{\dim \sigma}.$$
\\

Let $U_f$ be the affine hypersurface in $\gg_m^n$ defined by $f$. Define $\overline{U_f}:= \mathrm{Proj}(S_{\bigtriangleup}/(x_0f))$. It is the projective closure of $U_f$ in $\pp_{\bigtriangleup}$ defined over $\ff_q$. For any face $\sigma \subseteq \bigtriangleup$, we obtain an affine hypersurface $U_{f_{\sigma}}= \overline{U_f} \cap \gg_m^{\dim \sigma}$ in the algebraic torus $\gg_m^{\dim \sigma}$ defined by $f_{\sigma}= \sum_{j \in \sigma \cap J}a_jx^j$. Since $\overline{U_f}$ lies in $\pp_{\bigtriangleup}$, we still have the toric decomposition  
$$\overline{U_f}= \bigsqcup_{\emptyset \neq \sigma \subseteq \bigtriangleup} U_{f_{\sigma}}.$$
\\

Let $V$ be a variety over $\ff_q$, $V(\ff_{q^k})$ be the set of $\ff_{q^k}$-rational points of $V$ and $N_k(V)$ be its cardinality. The zeta function of $V$ is defined to be
$$Z(V/ \ff_q, t)= \exp \left( \sum_{k=1}^{\infty} N_k(V) \frac{t^k}{k} \right).$$ 
Let $\Psi : \ff_q \rightarrow \qq (\zeta_p)$ be a nontrivial additive character of $\ff_q$. Define the exponential sum of $x_0f$ to be 
$$S_k^*(x_0f)= \sum_{(x_0, x) \in (\ff_{q^k}^{*})^{n+1}} \Psi \circ \mathrm{Tr}_{\ff_{q^k} / \ff_q}(x_0f(x)) \in \qq (\zeta_p),$$ 
and the associated $L$-function of $x_0f$ to be 
$$L^*(x_0f, t)= \exp \left( \sum_{k=1}^{\infty} S_k^*(x_0f) \frac{t^k}{k} \right) \in \qq (\zeta_p)[[t]].$$
By theorem of Dwork-Bombieri-Grothendieck, $L^*(x_0f, t)$ is a rational function of $t$.\\

\noindent $\bf{Lemma \mbox{ 2.1.   }}$ Let $f \in \ff_{q^k}[x_1^{\pm 1}, \ldots x_n^{\pm 1}]$, $U_f$ the variety defined by $f=0$ in the torus $\gg_m^n$ and $N_k^*(U_f)$ the number of $\ff_{q^k}$-rational points of $U_f$. Then 
\begin{eqnarray*}
\sum_{(x_0, x) \in (\ff_{q^k}^{*})^{n+1}} \Psi \circ \mathrm{Tr}_{\ff_{q^k} / \ff_q} (x_0f(x))=q^kN_k^*(U_f)-N_k(\gg_m^n)
\end{eqnarray*}
where $x=(x_1, \ldots , x_n)$.\\

\noindent $\bf{Proof.}$ For any $x \in (\ff_{q^k}^*)^n$, we have 
$$\sum_{x_0 \in \ff_{q^k}} \Psi \circ \mathrm{Tr}_{\ff_{q^k} / \ff_q} (x_0f(x))= \left\{ \begin{array}{lcl}
q^k & \mbox{   if    } & f(x)=0 \mbox{   (i.e.,    } x \in U_f),\\
0 & \mbox{  } & \mbox{otherwise.} 
\end{array} \right.$$
This is a standard result from the theory of additive character sums. Thus 
$$\sum_{(x_0, x) \in \ff_{q^k} \times (\ff_{q^k}^{*})^n} \Psi \circ \mathrm{Tr}_{\ff_{q^k} / \ff_q} (x_0f(x))=q^kN_k^*(U_f).$$ 
In fact, $\sum_{(x_0, x) \in \ff_{q^k} \times (\ff_{q^k}^{*})^n} \Psi \circ \mathrm{Tr}_{\ff_{q^k} / \ff_q} (x_0f(x))$ is a counting function, counting how many $\ff_{q^k}$-rational points of $U_f$. Hence
\begin{eqnarray*} 
\sum_{(x_0, x) \in (\ff_{q^k}^{*})^{n+1}} \Psi \circ \mathrm{Tr}_{\ff_{q^k} / \ff_q} (x_0f(x)) &=& \sum_{(x_0, x) \in \ff_{q^k} \times (\ff_{q^k}^{*})^n}    \Psi \circ \mathrm{Tr}_{\ff_{q^k} / \ff_q} (x_0f(x))\\
& & - \sum_{(x_0, x) \in 0 \times (\ff_{q^k}^{*})^n} \Psi \circ \mathrm{Tr}_{\ff_{q^k} / \ff_q} (x_0f(x))\\
&=& q^kN_k^*(U_f)-(q^k-1)^n\\
&=& q^kN_k^*(U_f)-N_k(\gg_m^n).
\end{eqnarray*}
$\hfill$ $\Box$\\

Hence, we have 
$$Z(U_f, qt)=Z(\gg_m^n, t)L^*(x_0f, t).$$
Applying the above lemma to each face $\sigma \subseteq \bigtriangleup (f)$, we also have the analogous equations 
$$Z(U_{f_{\sigma}}, qt)=Z(\gg_m^{\dim \sigma}, t)L^*(x_0f_{\sigma}, t)$$
where $L^*(x_0f_{\sigma}, t)$ denotes the corresponding $L$-function with respect to $x_0f_{\sigma}$.\\

Since $\overline{U_f}= \bigsqcup_{\emptyset \neq \sigma \subseteq \bigtriangleup}U_{f_{\sigma}}$, $N_k(\overline{U_f})= \sum_{\emptyset \neq \sigma \subseteq \bigtriangleup} N_k(U_{f_{\sigma}})$. We have 
\begin{eqnarray*}
Z(\overline{U_f}, qt)= \prod_{\emptyset \neq \sigma \subseteq \bigtriangleup} Z(U_{f_{\sigma}}, qt)= \prod_{\emptyset \neq \sigma \subseteq \bigtriangleup} Z(\gg_m^{\dim \sigma}, t) \prod_{\emptyset \neq \sigma \subseteq \bigtriangleup} L^*(x_0f_{\sigma}, t).
\end{eqnarray*}
Thus it suffices to consider $\prod_{\emptyset \neq \sigma \subseteq \bigtriangleup} L^*(x_0f_{\sigma}, t)$.\\

Let $E(t)$ be the Artin-Hasse exponential series: 
$$E(t)= \exp \left( \sum_{i=0}^{\infty} \frac{t^{p^i}}{p^i} \right) \in (\zz_p \cap \qq)[[t]]$$
\\
\noindent Let $\gamma \in \Omega_1$ be a root of $\sum_{i=0}^{\infty} \frac{t^{p^i}}{p^i}=0$ satisfying ord$\gamma= \frac{1}{(p-1)}$, and consider 
$$\theta (t)=E(\gamma t)= \sum_{i=0}^{\infty} \lambda_i t^i \in \Omega_1[[t]].$$
The series $\theta (t)$ is a splitting function in Dwork's terminology \cite{Dwork}. In particular, its coefficients satisfy
\begin{eqnarray*}
\mbox{ord} \lambda_i \geq \frac{i}{p-1}, \mbox{      } \lambda_i \in \Omega_1.
\end{eqnarray*}
\\

Let ${\cal{O}}_0$ be the ring of integers of $\Omega_0$. Consider the following spaces of $p$-adic functions (where $b \in \rr$, $b \geq 0$, $c \in \rr$):
\begin{eqnarray*}
L(b, c) &=& \left \{ \sum_{(r, m) \in M(x_0f)} A_{r, m}x_0^rx^m|A_{r, m} \in \Omega_0, \mathrm{ord}A_{r, m} \geq br+c \right \} , \\
L(b) &=& \bigcup_{c \in \rr} L(b, c) ,\\
B({\cal{O}}_0) &=& \left \{ \sum_{(r, m) \in M(x_0f)} A_{r, m} \gamma^{r}x_0^rx^m|A_{r, m} \in {\cal{O}}_0, A_{r, m} \rightarrow 0 \mbox{   as   } r \rightarrow \infty \right \} , \\
B &=& \left \{ \sum_{(r, m) \in M(x_0f)} A_{r, m} \gamma^{r}x_0^rx^m| A_{r, m} \in \Omega_0, A_{r, m} \rightarrow 0 \mbox{   as   } r \rightarrow \infty \right \} .
\end{eqnarray*}

Observe that if $b>1/(p-1)$ then $L(b) \subseteq B \subseteq L(1/(p-1))$. If in addition $c \geq 0$, then $L(b, c) \subseteq B({\cal{O}}_0)$. Define a norm on $B$ as follows: If $\xi = \sum_{(r, m) \in M(x_0f)} A_{r, m} \gamma^{r} x_0^rx^m $   $\in B$, then define
$$\left\| \xi \right\| = \sup_{(r, m) \in M(x_0f)}|A_{r, m}|.$$
Note that $B({\cal{O}}_0)$ is a flat, separated, complete $\cal{O}$$_0$-module (\cite{Mon}, P.91).\\

Let $\widehat{f}= \sum_{j \in J} \widehat{a_j}x^j \in K[x_1^{\pm 1}, \ldots ,x_n^{\pm 1}]$ be the Teichm$\ddot{\mbox{u}}$ller lifting of $f$; i.e., $(\widehat{a_j})^q= \widehat{a_j}$ and the reduction of $\widehat{f}$ mod $p$ is $f$. Let 
$$F(f, x_0, x)= \prod_{j \in J} \theta \left ( \widehat{a_j}x_0x^j \right ) ,$$
$$F_a(f, x_0, x)= \prod_{i=0}^{a-1}F^{\tau^i}(f, x_0^{p^i}, x^{p^i})$$
where the map $\tau$ acts coefficient-wise on the power series $F(f, x_0, x)$. The estimate $\mathrm{ord}\lambda_i \geq i/(p-1)$ implies $F(f, x_0 ,x)$ and $F_a(f, x_0, x)$ are well-defined as formal Laurent series in $x_0, x_1^{\pm 1} , \ldots ,$   $ x_n^{\pm 1}$ with coefficients in $\Omega_0$; in fact 
$$F(f, x_0, x) \in L \left( \frac{1}{p-1}, 0 \right) , \mbox{             } F_a(f, x_0 ,x) \in L \left( \frac{p}{q(p-1)}, 0 \right) .$$
An easy calculation implies that $L(b)$ is a ring.\\

Define an operator $\psi$ on the formal Laurent series by 
$$\psi \left ( \sum_{(r, m) \in M(x_0f)}A_{r, m}x_0^rx^m \right ) = \sum_{(r, m) \in M(x_0f)}A_{pr, pm}x_0^rx^m.$$
where $pm$ denotes the $n$-tuple $(pm_1, \ldots , pm_n)$. This is the map on the power series which ignores about all $x_0^rx^m$-terms for which $p \nmid (r, m)$ and replaces $x_0^rx^m$ by $x_0^{r/p}x^{m/p}$ in the terms for which $p| (r, m)$. Here $p| (r, m)$ means that $p$ divides all of the entries in the integer vector $(r, m)$. One immediately deduces that $\psi(L(b, c)) \subseteq L(pb, c)$.\\

Let $\iota :L(p/(p-1)) \hookrightarrow L(p/q(p-1))$ be the canonical injection and denote by $\alpha$ the composition
$$L \left( \frac{p}{p-1} \right) \stackrel{\iota}{\hookrightarrow} L \left( \frac{p}{q(p-1)} \right) \stackrel{F_a(f, x_0, x)}{\longrightarrow} L \left( \frac{p}{q(p-1)} \right) \stackrel{\psi^a}{\rightarrow} L \left( \frac{p}{p-1} \right),$$
where the middle arrow means ``multiplication by $F_a(f, x_0, x)$''. It follows that the operator $\alpha = \psi^a \circ F_a(f, x_0, x)$ is an $\Omega_0$-linear endomorphism of the space $B$ and $L(b)$ for $0 \leq b \leq p/(p-1)$. Furthermore, the operator $\alpha_0= \tau^{-1} \circ \psi \circ F(f, x_0, x)$ is an $\Omega_1$-linear endomorphism of $B$ and $L(b)$ for $0 \leq b \leq p/(p-1)$ and is an $\Omega_0$-semilinear endomorphism of those spaces. It follows from Serre \cite{Serre} that the operators $\alpha^k$ and $\alpha_0^k$, acting on $B$ and $L(b)$ for $0<b \leq p/(p-1)$, have well-defined traces.  In addition, the Fredholm determinant $\det (I-t \alpha)$ and $\det (I-t \alpha_0)$ are $p$-adically entire (i.e., convergent for all $t \in \Omega$). Dwork trace formula asserts that 
\begin{eqnarray}
S_k^*(x_0f)=(q^k-1)^{n+1} \mathrm{Tr}(\alpha^k),
\end{eqnarray}
where $\alpha$ acts either on $B$ or on some $L(b)$, $0<b \leq p/(p-1)$. The non-trivial additive character $\Psi$ implicit on the left-hand side of (2) is given by 
$$\Psi (t)= \theta (1)^{ \mathrm{Tr}_{\ff_q/ \ff_p}(t)}$$ 
for $t \in \ff_q$ (by \cite{Dwork}, Lemma 4.1, $\theta (1)$ is a primitive $p$-th root of unity). Equivalently, one can define an operator $\delta$ on formal power series with constant term $1$ by $g(t)^{\delta}=g(t)/g(qt)$, and then (2) takes the form
\begin{eqnarray}
L^*(x_0f, t)^{(-1)^n}= \det (I-t \alpha)^{\delta^{n+1}}
\end{eqnarray}
via the relationship $\det (I-t \alpha)= \exp (- \sum_{k=1}^{\infty} \mathrm{Tr}(\alpha^k)t^k/k)$.\\

\section{Dwork Cohomology}
\

Let $f$ be a Laurent polynomial. Define the new polynomials $F(x_0, x)=x_0f(x)-1$ and 
$$F_i(x_0, x)=x_i \frac{\partial}{\partial x_i}F(x_0, x) \mbox{   for   } 0 \leq i \leq n.$$
\\

\noindent $\bf{Definition.}$  A Laurent polynomial $f$ is called $\bigtriangleup$-regular if $\bigtriangleup (f)= \bigtriangleup$ and for every $l$-dimensional face $\sigma \subseteq \bigtriangleup$   $(l>0)$ the polynomial equations
\begin{eqnarray*}
f_{\sigma}(x)=x_1 \frac{ \partial f_{\sigma}}{ \partial x_1} = \ldots =x_n  \frac{ \partial f_{\sigma} } { \partial x_n} =0
\end{eqnarray*}
have no common solution in $\gg_m^n$.\\

For $0 \leq l \leq n$, the sum
$$\bigoplus_{|T|=l} S_{\bigtriangleup}e_T$$ 
runs over all subsets $T \subseteq \{ 0, \ldots , n \}$ of cardinality $l$. If $T= \{ i_1, \ldots ,i_l \}$ with $i_1 < \ldots < i_l$, then $e_T$ is the abbreviated notation for $e_{i_1} \wedge \ldots \wedge e_{i_l}$ (the $\{ e_i \}_{0 \leq i \leq n}$ being a set of formal symbols). 
The boundary map $\overline{\partial_l}: \bigoplus_{|T|=l} S_{\bigtriangleup}e_T \rightarrow \bigoplus_{|T|=l-1} S_{\bigtriangleup}e_T$ is given by:   
$$\overline{\partial_l}(\xi e_T)= \sum_{j=1}^{l} (-1)^{j-1} F_j \xi e_{T- \{ i_j \}}$$
where $\xi \in S_{\bigtriangleup}$ and $T= \{ i_1, \ldots ,i_l \}$ with $i_1 < \ldots < i_l$. We have the Koszul complex $K_*(f)$ on the elements $F_0, \ldots , F_n \in S_{\bigtriangleup}$:
$$K_*(f):0 \rightarrow \bigoplus_{|T|=n+1} S_{\bigtriangleup}e_T \stackrel{\overline{\partial_{n+1}}}{\rightarrow} \ldots \rightarrow \bigoplus_{|T|=1} S_{\bigtriangleup}e_T \stackrel{\overline{\partial_1}}{\rightarrow} S_{\bigtriangleup} \rightarrow 0.$$
\\

\noindent $\bf{Theorem \mbox{   3.1.   } (\cite{Bat}}, \mbox{  Theorem 4.8  })$ Let $f$ be a Laurent polynomial and $F=x_0f(x)-1$. Then the following conditions are equivalent:\\
(i) $f$ is $\bigtriangleup$-regular;\\
(ii) $F$ is $\bigtriangleup$-regular;\\
(iii) the elements $\{ F_0, F_1, \ldots , F_n \}$ give rise to a regular sequence in $S_{\bigtriangleup}$;\\
(iv) the homology groups $H_i(f)$ of the Koszul complex $K_*(f)$ are zero for positive $i$;\\
(v) the dimension of $H_0(f)$ is $(n+1)! \mathrm{Vol} (x_0f)$.\\  

Let $\widehat{F} (x_0, x)$ be the Teichm$\ddot{\mbox{u}}$ller lifting of $F(x_0, x)$. Define $\widehat{F_i} =x_i \partial \widehat{F} / \partial x_i$ and $\gamma_l= \sum_{i=0}^l \gamma^{p^i} /p^i$, which by definition of $\gamma$ satisfies 
$$\mathrm{ord}_p \gamma_l \geq \frac{p^{l+1}}{p-1}-l-1.$$
For $i=0, \ldots , n$, define differential operators $\widehat{D_i}$ by 
$$\widehat{D_i}=E_i+ \widehat{H_i},$$
where $E_i=x_i \partial / \partial x_i$ and 
$$\widehat{H_i}(x_0, x)= \sum_{l=0}^{\infty} \gamma_lp^l \widehat{F_i}^{\tau^l}(x_0^{p^l}, x^{p^l}) \in L \left( \frac{p}{p-1}, -1 \right) .$$
The $\widehat{H_i}$ and $\widehat{D_i}$ operate on $B$ and on $L(b)$ for $b \leq p/(p-1)$. One verifies that they commute with one another. Furthermore, by \cite{Dwork} equation (4.35), 
\begin{eqnarray}
\alpha \circ \widehat{D_i} =q \widehat{D_i} \circ \alpha
\end{eqnarray}
for $i=0, \ldots , n$ as operators on $B$ or on $L(b)$, $0<b \leq p/(p-1)$.\\

Let $K_*(B, \{ \widehat{D_i} \}_{i=0}^n)$ be the Koszul complex on $B$ formed from $\widehat{D_0}, \ldots , \widehat{D_n}$, i.e., for $0 \leq l \leq n$, 
$$K_l(B)= \bigoplus_{|T|=l}Be_T,$$
where the sum runs over all subsets $T \subseteq \{ 0, \ldots , n \}$ of cardinality $l$. The boundary map $\partial_l:K_l \rightarrow K_{l-1}$ is given by: 
$$\partial_l(\xi e_T)= \sum_{j=1}^{l} (-1)^{j-1} \widehat{D_{i_j}}(\xi)e_{T- \{ i_j \}}$$
where $\xi \in B$ and $T= \{ i_1, \ldots ,i_l \}$ with $i_1 < \ldots < i_l$. Define an endomorphism $\alpha_l:K_l \rightarrow K_l$ by 
$$\alpha_l= \bigoplus_{|T|=l}q^l \alpha.$$
\\

\noindent Equation (4) implies that $\alpha_*$ is a chain map on the Koszul complex $K_*$. 
$$\begin{CD}
0 @> >> K_{n+1} @> \partial_{n+1} >> K_n @> \partial_n >> \ldots @> \partial_1 >> K_0 @> >>0 \\
@V VV @V \alpha_{n+1} VV @V \alpha_n VV @V \ldots VV @V \alpha_0 VV @V VV\\
0 @> >> K_{n+1} @> \partial_{n+1} >> K_n @> \partial_n >> \ldots @> \partial_1 >> K_0 @> >>0 \\
\end{CD}$$

\noindent It then follows from (3) that
\begin{eqnarray*} 
L^*(x_0f, t)^{(-1)^n}= \prod_{l=0}^{n+1} \det (I-t \alpha_l|K_l)^{(-1)^l}.
\end{eqnarray*}
Passing to homology, we conclude that
\begin{eqnarray*} 
L^*(x_0f, t)^{(-1)^n}= \prod_{l=0}^{n+1} \det (I-t H_l(\alpha_*)|H_l(K_*))^{(-1)^l},
\end{eqnarray*}
where $H_l(\alpha_*)$ denotes the endomorphism of $H_l(K_*)$ induced by $\alpha_l$.\\

Some of the main results of \cite{AS1} are summarized in the following theorem.\\

\noindent $\bf{Theorem \mbox{ 3.2.   }}$ Let $f \in \ff_q[x_1^{\pm 1}, \ldots , x_n^{\pm 1}]$ be $\bigtriangleup$-regular and $\dim \bigtriangleup (f)=n$. The complex $K_*(B, \{ \widehat{D_i} \}_{i=0}^n)$ is acyclic in positive dimension. In addition, $H_0(K_*(B, \{ \widehat{D_i} \}_{i=0}^n))$ has dimension $(n+1)! \mathrm{Vol}(x_0f)$ and $H_0(\alpha_*)$ is invertible. Hence 
\begin{eqnarray*}
L^*(x_0f, t)^{(-1)^n}= \det (I-tH_0(\alpha_*)|H_0(K_*(B, \{ \widehat{D_i} \}_{i=0}^n))) \in \zz[t]
\end{eqnarray*}
is a polynomial of degree $(n+1)! \mathrm{Vol}(\bigtriangleup(x_0f))$.\\

\noindent $\bf{Proof.}$ All statements follows immediately from Theorem 2.9, Theorem 2.18, Corollary 2.19, and Theorem 3.13 of \cite{AS1}. $\hfill$ $\Box$\\

\section{Simplex $\bigtriangleup (f)$}
\   
   
Suppose $H$ is a hyperplane in $\rr^{n+1}$ passing through the origin: $\sum_{j=0}^na_jx_j=0$. Suppose in addition the $a_j$ are rational. We denote by $\widehat{D_H}$ the differential operator on $B$ defined by $\widehat{D_H}= \sum_{j=0}^na_j \widehat{D_j}$.\\

Suppose that $\bigtriangleup(f)$ is an $n$-dimensional simplex with vertices $\{ v_0, \ldots ,v_n \} \subseteq \rr^n$. Then so is $\bigtriangleup(x_0f)$ with vertices $ \{ (0, \ldots ,0), (1, v_0), \ldots ,(1, v_n) \} \subseteq \rr^{n+1}$. Let $H_i$ be a hyperplane through the codimension-one faces of $\bigtriangleup(x_0f)$ that contain $\{ (0, \ldots ,0),$ $ (1, v_0), \ldots , \widehat{(1, v_i)} , \ldots ,(1, v_n) \}$. Then the normal vectors of $H_0, \ldots ,H_n$ form a basis in $\rr^{n+1}$. We may assume the equation of each $H_i$ has the form
\begin{eqnarray}
\sum_{j=0}^na_{ij}x_j=0,     
\end{eqnarray}
where $a_{ij} \in \zz$ for all $i$, $j$ and $\{ a_{ij} \}_{j=0}^n$ have greatest common divisor $1$. The differential operators $\{ \widehat{D_{H_i}} \}_{i=0}^n$ on $B$ commute with one another and $\alpha \circ \widehat{D_{H_i}}=q \widehat{D_{H_i}} \circ \alpha$.\\

\noindent $\bf{Lemma \mbox{ 4.1.   }}$ Let $f \in \ff_q[x_1^{\pm 1}, \ldots , x_n^{\pm 1}]$ be $\bigtriangleup$-regular and $F(x_0, x)=x_0f(x)-1$. Suppose $p \nmid \det(a_{ij})$ where the matrix $(a_{ij})$ is defined as in (5). Then Theorem 3.2 is valid with $\{ \widehat{D_i} \}_{i=0}^n$ replaced by $\{ \widehat{D_{H_i}} \}_{i=0}^n$.\\

\noindent $\bf{Proof.}$ Since for each face $\tau$ of $\bigtriangleup (x_0f)$,
$$\left( \sum_{j=0}^na_{ij}x_j \partial F/ \partial x_j \right)_{\tau}= \sum_{j=0}^na_{ij}x_j \partial F_{\tau}/ \partial x_j.$$
The hypothesis that $f$ is $\bigtriangleup$-regular and that $p \nmid \det(a_{ij})$ implies that for each face $\tau$ of $\bigtriangleup(x_0f)$, the polynomials $(\sum_{j=0}^na_{ij}x_j \partial F/ \partial x_j  )_{\tau}$, $i=0, \ldots , n$, have no common zero in $(\ff_q^*)^{n+1}$. One can repeat the arguments of \cite{AS1} with $x_i \partial F/ \partial x_i$ replaced by $\sum_{j=0}^na_{ij}x_j \partial F/ \partial x_j$. $\hfill$ $\Box$\\  

For $j=0, \ldots , n$, we define $\Omega_0$-endomorphism $\theta_j:B \rightarrow B$ by 
\begin{eqnarray}
\theta_j \left( \sum_{(r, m) \in M(x_0f)} A_{r, m} \gamma^r x_0^rx^m \right) = \sum_{(r, m) \in M(x_0f) \cap H_j} A_{r, m} \gamma^r x_0^rx^m.
\end{eqnarray}
Note that $\theta_j$ is also a ring homomorphism of $B$ to itself.\\

\noindent $\bf{Lemma \mbox{ 4.2.   }}$ With the notation as above, we have 
\begin{eqnarray} 
\widehat{D_{H_j}}(B) \subseteq \ker (\theta_j)  &\mbox{   for   }& 0 \leq j \leq n\\
\theta_j^2 = \theta_j &\mbox{   for   } & 0 \leq j \leq n\\
\theta_j \circ \widehat{D_{H_i}} \circ \theta_j = \theta_j \circ \widehat{D_{H_i}} &\mbox{   for   }& 0 \leq i, j \leq n.
\end{eqnarray}

\noindent $\bf{Proof.}$ Let $\sum_{i=0}^n a_{ji}x_i=0$ be the defining equation of $H_j$. Write $x^u=x_0^{u_0} \cdots x_n^{u_n}$ where $u \in M(x_0f)$. Then $\widehat{D_{H_j}}= \sum_{i=0}^n a_{ji}(E_i+ \widehat{H_i}(x_0, x))$. Consider $\theta_j \circ (\sum_{i=0}^n a_{ji}E_i)$, we have 
\begin{eqnarray*}
\sum_{i=0}^n a_{ji}E_i(x^u)= \left( \sum_{i=0}^n a_{ji}u_i \right) x^u &=& 0 \mbox{  if   } u \in H_j \mbox{  and   }\\
\theta_j \left( \sum_{i=0}^n a_{ji}E_i(x^u) \right) = \left( \sum_{i=0}^n a_{ji}u_i \right) \theta_j(x^u) &=& 0 \mbox{  if   } u \notin H_j.
\end{eqnarray*}
Hence $\theta_j \circ (\sum_{i=0}^n a_{ji}E_i)(x^u)=0$ for any $u \in M(x_0f)$.\\

In order to show $\theta_j \circ (\sum_{i=0}^n a_{ji} \widehat{H_i}(x_0, x)) \cdot (x^u)=0$, it suffices to show $\theta_j \circ (\sum_{i=0}^n a_{ji} (x_i \frac {\partial} {\partial x_i} x^v)) \cdot (x^u)=0$ where $x^v=x_0^{v_0} \cdots x_n^{v_n}$ is a monomial of $F$. Clearly, 
$$\theta_j \circ \left( \sum_{i=0}^n a_{ji} (x_i \frac {\partial} {\partial x_i} x^v) \cdot x^u \right) = \theta_j \circ \left( \sum_{i=0}^n a_{ji}v_i  x^v \cdot x^u \right) = \left( \sum_{i=0}^n a_{ji}v_i \right) \theta_j (x^{u+v}).$$ 
If $v \in H_j$, then $\sum_{i=0}^n a_{ji}v_i=0$. If $v \notin H_j$, then $u+v \notin H_j$. Hence $\theta_j(x^{u+v})=0$. (7) is proved. (8) is clear.\\

The equation (9) is satisfied by $x^u$ where $u \in H_j$. For $u \notin H_j$, $\theta_j \circ \widehat{D_{H_i}} \circ \theta_j=0$. Let $\sum_{j=0}^na_{ij}x_j=0$ be the defining equation of $H_i$. Since $(\sum_{j=0}^na_{ij}E_j)(x^u)$ is a scalar multiple of $x^u$, $\theta_j (\sum_{j=0}^na_{ij}E_j)(x^u)=0$. Clearly, the exponents of the monomials of $\widehat{H_i}(x_0, x) \cdot x^u \notin H_j$  for all $i$ as $u \notin H_j$. That means $\theta_j \widehat{H_i}(x_0, x)(x^u)=0$ for all $i$. Hence $\theta_j \circ \widehat{D_{H_i}}=0$. $\hfill$ $\Box$\\

Let $S= \{ 0, 1, \ldots , n \}$. For $A \subseteq S$, let $\theta_A$ be the composition of all $\theta_j$ for $j \in A$. Let $B_A= \theta_A(B)$. The maps $\widehat{D_{H_j}^A} \stackrel{def}{:=} \theta_A \circ \widehat{D_{H_j}}$ for $j \in S$ are stable on $B_A$. Thus we may form the Koszul complex $K_*(B_A, \{ \widehat{D_{H_j}^A} \}_{j \in S-A})$ which we abbreviate by $K_*(A)$. Note that we ignore the differential operators $\widehat{D_{H_j}^A}$ for $j \in A$ since by (7) they are identically zero on $B_A$. For each subset $A^{'} \subseteq S-A$, we define $B_{A, A^{'}}= \bigcap_{j \in A^{'}}( \mathrm{ker} \theta_j|B_A)$ and a subcomplex $K_*(A, A^{'})$ of $K_*(A)$ by setting for $0 \leq l \leq n+1-|A|$
\begin{eqnarray}
K_l(A, A^{'})= \bigoplus_{T \subseteq S-A, |T|=l} \left( \bigcap_{j \in A^{'}-T}(ker \theta_j|B_A) \right)e_T,
\end{eqnarray}
where the sum runs over all subsets $T \subseteq S-A$ of cardinality $l$, and if $A^{'}-T= \emptyset$ we take $``\bigcap_{j \in A^{'}-T}( \mathrm{ker} \theta_j|B_A)"$ to mean $``B_A"$. It is straightforward to check that the boundary maps are stable on the submodules, hence they define a subcomplex. Note that $K_*(A)=K_*(A, \emptyset)$ and that $K_*( \emptyset)$ is the complex $K_*(B,  \{ \widehat{D_{H_j}} \}_{j=0}^n)$.\\

We can define an action of Frobenius on these complexes. Since $\theta_j$ commutes with $\psi$ and with multiplication by $F_a(f, x_0, x)$, we may define an endomorphism 
$$\alpha_A= \psi^a \circ \theta_A(F_a(f, x_0, x)):B_A \rightarrow B_A$$
that satisfies
$$\alpha_A \circ D_{H_i}^A=qD_{H_i}^A \circ \alpha_A$$
for $i=0,1, \ldots , n$. As before, this induces an endomorphism $\alpha_A:K_*(A) \rightarrow K_*(A)$ by setting
$$\alpha_{A,l}= \bigoplus_{|T|=l}q^l \alpha_A:K_l(A) \rightarrow K_l(A).$$
It is straightforward to check that $\alpha_A$ induces endomorphism $\alpha_{A, A^{'}}$ of the subcomplexes $K_*(A, A^{'})$ for all $A^{'} \subseteq S-A$, as well. We denote by $H_*(\alpha_A)$ and $H_*(\alpha_{A, A^{'}})$ the endomorphisms of the associated homology groups.\\

We now explain the arithmetic significance of these complexes. For $f= \sum_{j \in J} a_jx^j $     $\in \ff_q[x_1^{\pm 1}, \ldots , x_n^{\pm 1}]$ and $A \subseteq S$, define 
$$f_A= \sum_{j \in J, (1, j) \in \bigcap_{i \in A} H_i}a_jx^j \in \ff_q[x_1^{\pm 1}, \ldots , x_n^{\pm 1}].$$ 
We denote the corresponding exponential sums on $\gg_m^{n+1-|A|}$ by $S_k^*(x_0f_A)$ and the corresponding $L$-function by $L^*(x_0f_A, t)$.\\

In Section 3 we associate to $f$ a Koszul complex $K_*(B)$ with Frobenius operator $\alpha$ satisfying Theorem 3.2 when $f$ is $\bigtriangleup$-regular. The same construction can be applied to each $f_A$, and we denote the resulting Koszul complex and Frobenius operator by $K_*(B(x_0f_A))$ and $\alpha_{x_0f_A}$, respectively. Let $\widehat{D_{A,i}}$, $i=0, \ldots ,n$ be the corresponding differential operators. When $f$ is $\bigtriangleup$-regular, so are all the $f_A$, hence Theorem 3.2 holds for Koszul complex $K_*(B(x_0f_A), \{ \widehat{D_{A,i}} \}_{i=0}^n)$. In addition, if $p \nmid \det (a_{ij})$ the proof of Lemma 4.1 show that Theorem 3.2 holds for the Koszul complex $K_*(B(x_0f_A), \{ \widehat{D_{A, H_i}} \}_{i=0}^n)$, where 
$$\widehat{D_{A, H_i}} \stackrel{def}{:=} \sum_{j=0}^na_{ij} \widehat{D_{A,j}}.$$

Comparing $K_*(B(x_0f_A), \{ \widehat{D_{A,i}} \}_{i=0}^n)$ with the complex $K_*(A)$ constructed in this section, one sees that $K_l(B(x_0f_A))$ and $K_l(A)$ are naturally identified as Banach spaces with action of Frobenius. Under this identification of spaces, the differential operator $\widehat{D_{A, H_i}}$ is identified with $\widehat{D_{H_i}^A}$ for $i=0 , \ldots , n$. Thus $K_*(B(x_0f_A), \{ \widehat{D_{A, H_i}} \}_{i=0}^n)$ and $K_*(B_A, \{ \widehat{D_{H_i}^A} \}_{i=0}^n)$ are isomorphic as complexes with action of Frobenius. The hypothesis that $\bigtriangleup(f)$ be simplex implies that $\dim \bigtriangleup (x_0f_A)=n+1-|A|$ for each $A \subseteq S$. Keeping in mind that $\widehat{D_{H_i}^A}=0$ (as operator on $B_A$) for $i \in A$, Theorem 3.2 applied to $x_0f_A$ gives:\\   

\noindent $\bf{Lemma \mbox{ 4.3.   }}$ Let $f \in \ff_q[x_1^{\pm 1}, \ldots , x_n^{\pm 1}]$ be $\bigtriangleup$-regular and suppose $\bigtriangleup (f)$ is an $n$-dimensional simplex and $p \nmid \det (a_{ij})$ where the matrix is defined as in (5). For all $A \subseteq S$, the complex $K_*(A)$ is acyclic in dimension $>0$, $\dim H_0(K_*(A))=(n+1-|A|)! \mathrm{Vol}(x_0f_A)$, and $H_0(\alpha_A)$ is invertible. Define
$$P_A(t)= \det (I-t \alpha_A|B_A).$$
Then $\det(I-tH_0(\alpha_A)|H_0(K_*(A)))$ is a polynomial in $\zz[t]$ of degree $(n+1-|A|)! $              $\mathrm{Vol}(x_0f_A)$ and 
\begin{eqnarray*}
L^*(x_0f_A, t)^{(-1)^{n-|A|}}=P_A(t)^{\delta^{n+1-|A|}}.
\end{eqnarray*}
\\

Let $A \subseteq S$, $A^{'} \subseteq S-A$, and $j \in A^{'}$. By (8), the map $\theta_j:B_{A, A^{'}- \{ j \}} \rightarrow B_{A \cup \{ j \}, A^{'}- \{ j \}}$ is surjective. The kernel of $\theta_j$ is naturally identified with $B_{A, A^{'}}$. Hence there is a short exact sequence  
\begin{eqnarray}
0 \rightarrow B_{A, A^{'}} \rightarrow B_{A, A^{'}- \{ j \}} \stackrel{\theta_j}{\rightarrow} B_{A \cup \{ j \}, A^{'}- \{ j \}} \rightarrow 0.
\end{eqnarray}
The surjective map $\theta_j:B_{A, A^{'}- \{ j \}} \rightarrow B_{A \cup \{ j \}, A^{'}- \{ j \}}$ induces a surjective map $\theta_{j, *}$ of complexes  
\begin{eqnarray*}
\theta_{j, *}:K_*(A, A^{'}- \{ j \}) \rightarrow K_*(A \cup \{ j \}, A^{'}- \{ j \})
\end{eqnarray*}
with kernel $K_*(A, A^{'})$.  Hence there is a short exact sequence of complexes with Frobenius action 
\begin{eqnarray}
0 \rightarrow K_*(A, A^{'}) \rightarrow K_*(A, A^{'}- \{ j \}) \stackrel{\theta_{j, *}}{\rightarrow} K_*(A \cup \{ j \}, A^{'}- \{ j \}) \rightarrow 0.
\end{eqnarray}
\\
   
Let $P_{A, A^{'}}(t)= \det (I-t \alpha_{A, A^{'}}|B_{A, A^{'}})$. Note that $P_{A, \emptyset}(t)=P_A(t)$. The short exact sequence (11) gives
\begin{eqnarray}
P_{A, A^{'}}(t)= \frac{P_{A, A^{'}- \{ j \}}(t)}{P_{A \cup \{ j \},A^{'}- \{ j \}} (t)}.
\end{eqnarray}
\\

We summarize these results:\\
 
\noindent $\bf{Lemma \mbox{ 4.4.   }}$ Let $f \in \ff_q[x_1^{\pm 1}, \ldots , x_n^{\pm 1}]$ be $\bigtriangleup$-regular and suppose $\bigtriangleup (f)$ is an $n$-dimensional simplex and $p \nmid \det (a_{ij})$ where the matrix $(a_{ij})$ is defined as in (5). For all $A \subseteq S$ and $A^{'} \subseteq S-A$, we have    
\begin{eqnarray*}
P_{A, A^{'}}(t)^{(-1)^{|A|}}= \prod_{A \subseteq C \subseteq A \cup A^{'}} P_C(t)^{(-1)^{|C|}}.
\end{eqnarray*}
\\

\noindent $\bf{Proof.}$ Use induction on $|A^{'}|$. For $|A^{'}|=0$, that means $A^{'}= \emptyset$ and $C=A$. 
$$P_{A, \emptyset}(t)^{(-1)^{|A|}}=P_A(t)^{(-1)^{|A|}}= \prod_{A \subseteq C \subseteq A \cup \emptyset} P_C(t)^{(-1)^{|C|}}.$$
Suppose $A^{'}$ is non-empty. For $j \in A^{'}$ by (13), we have
\begin{eqnarray*}
P_{A, A^{'}}(t)^{(-1)^{|A|}} &=& \frac{P_{A, A^{'}- \{ j \}}(t)^{(-1)^{|A|}}}{P_{A \cup \{ j \},A^{'}- \{ j \}}(t)^{(-1)^{|A|}}}\\
&=& \left( P_{A, A^{'}- \{ j \}}(t)^{(-1)^{|A|}} \right) \left( P_{A \cup \{ j \},A^{'}- \{ j \}}(t)^{(-1)^{|A|+1}} \right)\\
& & \mbox{   by induction hypothesis}\\
&=& \left( \prod_{A \subseteq C_1 \subseteq A \cup A^{'}- \{ j \}} P_{C_1}(t)^{(-1)^{|C_1|}} \right) \left( \prod_{A \cup \{ j \} \subseteq C_2 \subseteq A \cup A^{'}} P_{C_2}(t)^{(-1)^{|C_2|}} \right)\\ 
&=& \left( \prod_{A \subseteq C \subseteq A \cup A^{'}, j \notin C}  P_C(t)^{(-1)^{|C|}} \right) \left( \prod_{A \subseteq C \subseteq A \cup A^{'}, j \in C} P_C(t)^{(-1)^{|C|}} \right)\\
&=& \prod_{A \subseteq C \subseteq A \cup A^{'}} P_C(t)^{(-1)^{|C|}}.
\end{eqnarray*}
$\hfill$ $\Box$\\ 

\noindent Hence we have the following lemma.\\
   
\noindent $\bf{Lemma \mbox{ 4.5.   }}$ Let $f \in \ff_q[x_1^{\pm 1}, \ldots , x_n^{\pm 1}]$ be $\bigtriangleup$-regular and suppose $\bigtriangleup (f)$ is an $n$-dimensional simplex and $p \nmid \det (a_{ij})$ where the matrix $(a_{ij})$ is defined as in (5). For all $A \subseteq S$ and $A^{'} \subseteq S-A$, we have    
\begin{eqnarray*}
P_A(t)= \prod_{S-A^{'} \supseteq A} P_{S-A^{'}, A^{'}}(t)
\end{eqnarray*}
\\

\noindent $\bf{Proof.}$ By Lemma 4.4 we have 
\begin{eqnarray*} 
\prod_{S-A^{'} \supseteq A}P_{S-A^{'}, A^{'}}(t) &=& \prod_{S-A^{'} \supseteq A} \prod_{S-A^{'} \subseteq C \subseteq S} P_C(t)^{(-1)^{|C|+|S-A^{'}|}}\\
&=& \prod_{A \subseteq C} \prod_{A \subseteq S-A^{'} \subseteq C} P_C(t)^{(-1)^{|C|+|S-A^{'}|}}\\
&=& \prod_{A \subseteq C} P_C(t)^{(-1)^{|C|} \left( \sum_{A \subseteq S-A^{'} \subseteq C} (-1)^{|S-A^{'}|} \right) }\\
&=& P_A(t)
\end{eqnarray*}
since 
$$\sum_{A \subseteq S-A^{'} \subseteq C} (-1)^{|S-A^{'}|}= \left\{ \begin{array}{lcl}
(-1)^{|A|} & \mbox{   if    } & C=A,\\
0 & \mbox{   if   } & C \neq A.
\end{array} \right.$$
$\hfill$ $\Box$\\

The geometric meaning of above can be interpreted as follows: $B_A$ consists of monomials $x_0^rx^m$ where $(r, m)$ lies in the cone $M(x_0f_A)= \bigcap_{j \in A} H_j$. For $A^{'} \subseteq S-A$, $B_{A, A^{'}}$ consists of monomials $x_0^rx^m$ where $(r, m)$ lies in the cone $M(x_0f_A)$ but not in $H_j$, $j \in A^{'}$. In particular, $B_{A, S-A}$ consists of monomials $x_0^rx^m$ where $(r, m)$ lies in the interior of the cone $M(x_0f_A)$. In the cone $M(x_0f_A)$, the map $\theta_j:B_{A, A^{'}- \{ j \}} \rightarrow B_{A \cup \{ j \}, A^{'}- \{ j \}}$ in (6) is a projection onto $H_j$ which sends the monomials $x_0^rx^m$ to $0$ where $(r, m) \notin H_j$. That is why such monomials $x_0^rx^m \in B_{A, A^{'}- \{ j \}}$, with $(r,m) \notin H_j$, precisely lie in $B_{A, A^{'}}$ and have the exact sequence (11). Deleting the points $(r, m)$ not in $H_j$ where $x_0^rx^m \in B_{A, A^{'}- \{ j \}}$ gives (13). Lemma 4.4 is followed by the inclusion and exclusion principle and Lemma 4.5 by the boundary decomposition theorem \cite{Wan}. \\

\noindent $\bf{Lemma \mbox{ 4.6.   }}$ Suppose that for all $A \subseteq S$,
$$H_i(K_*(A, \emptyset))=0 \mbox{   for   } i>0.$$
Then for all $A^{'} \subseteq S-A$,
$$H_i(K_*(A, A^{'}))=0 \mbox{   for   } i>0.$$
In particular,
$$H_i(K_*(\emptyset , S))=0 \mbox{   for    } i>0.$$
Furthemore, there is an injection
$$H_0(K_*(\emptyset ,S)) \hookrightarrow H_0(K_*(\emptyset , \emptyset)).$$
\\

\noindent $\bf{Proof.}$ The proof is by induction on $|A^{'}|$, the case $|A^{'}|=0$ being the hypothesis of the theorem. From (12) we get the long exact sequence 
$$\cdots \rightarrow H_{i+1}(K_*(A \cup \{ j \} , A^{'}- \{ j \})) \rightarrow H_i(K_*(A , A^{'})) \rightarrow H_i(K_*(A , A^{'}- \{ j \})) \rightarrow \cdots $$
where $j \in A^{'}$. By induction hypothesis, the two ``outside'' homology groups vanish for $i>0$. Therefore $H_i(K_*(A, A^{'}))=0$ for $i>0$ also. The exact homology sequence associated to (12) then reduces to the short exact sequence 
\begin{eqnarray}
0 \rightarrow H_0(K_*(A, A^{'})) \rightarrow H_0(K_*(A, A^{'}- \{ j \})) \rightarrow H_0(K_*(A \cup \{ j \} , A^{'}- \{ j \})) \rightarrow 0.
\end{eqnarray}
In particular, there is an injection $H_0(K_*(A, A^{'})) \hookrightarrow H_0(K_*(A, A^{'}- \{ j \})).$ The existence of the injection asserted by the lemma then follows by induction on $|A^{'}|$. $\hfill$ $\Box$\\


\noindent $\bf{Theorem \mbox{ 4.7.   }}$ Let $f \in \ff_q[x_1^{\pm 1}, \ldots , x_n^{\pm 1}]$ be $\bigtriangleup$-regular and $\overline{U_f}$ the toric projective closure defined by $f$. Suppose $\bigtriangleup (f)$ is an $n$-dimensional simplex and $p \nmid \det (a_{ij})$ where the matrix $(a_{ij})$ is defined as in (5). Then
\begin{eqnarray*}
Z(\overline{U_f}, qt) = \frac{\det \left( I-t H_0(\alpha_{\emptyset, S})|H_0(K_*(\emptyset, S)) \right)^{(-1)^n}}{(1-qt) \cdots (1-q^nt)} 
\end{eqnarray*}
\\

\noindent $\bf{Proof.}$ By Lemma 4.3 and Lemma 4.5, we have   
\begin{eqnarray*} 
\prod_{A \subseteq S} L^*(x_0f_A, t)^{(-1)^n} &=&  \left( \prod_{A \subseteq S} L^*(x_0f_A, t)^{(-1)^{n-|A|}} \right)^{(-1)^{|A|}}\\
&=& \prod_{A \subseteq S} \left( P_A(t)^{\delta^{n+1-|A|}} \right)^{(-1)^{|A|}}\\ 
&=& \prod_{A \subseteq S} \prod_{i=0}^{n+1-|A|} P_A(q^it)^{(-1)^{i+|A|} \binom{n+1-|A|}{i}}\\
&=& \prod_{A \subseteq S} \prod_{i=0}^{n+1-|A|} \prod_{S-A^{'} \supseteq A} P_{S-A^{'}, A^{'}}(q^it)^{(-1)^{i+|A|} \binom{n+1-|A|}{i}}\\
&=& \prod_{A \subseteq S} \prod_{S-A^{''} \supseteq A} \prod_{S-A^{'} \supseteq A} P_{S-A^{'}, A^{'}}(q^{|A^{''}|}t)^{(-1)^{|A|+|A^{''}|}}\\
&=& \prod_{A^{''}, A^{'}} P_{S-A^{'}, A^{'}}(q^{|A^{''}|}t)^{(-1)^{|A^{''}|} \left( \sum_{A \subseteq (S-A^{'}) \cap (S-A^{''})} (-1)^{|A|} \right) }\\
\end{eqnarray*}
where $A$ corresponds to the face $\sigma = \bigcap_{i \in A} H_i$ of $\bigtriangleup (x_0f)$, while $S-A^{'} \supseteq A$ and $S-A^{''} \supseteq A$ correspond to the subfaces $\sigma^{'} = \bigcap_{i \in S-A^{'}} H_i$ and $\sigma^{''} = \bigcap_{i \in S-A^{''}} H_i$ of $\sigma$. Since   
$$\sum_{A \subseteq (S-A^{'}) \cap (S-A^{''})} (-1)^{|A|} = \left\{ \begin{array}{lcl}
1 & \mbox{   if    } & (S-A^{'}) \cap (S-A^{''}) = \emptyset,\\
0 & \mbox{   if    } & (S-A^{'}) \cap (S-A^{''}) \neq \emptyset
\end{array} \right.$$
\begin{eqnarray*}
\prod_{A \subseteq S}L^*(x_0f_A, t)^{(-1)^n} &=& \prod_{(S-A^{'}) \cap (S-A^{''})= \emptyset} P_{S-A^{'}, A^{'}}(q^{|A^{''}|}t)^{(-1)^{|A^{''}|} }\\
&=& \prod_{A^{''}} \prod _{S-A^{'} \subseteq A^{''}} P_{S-A^{'}, A^{'}}(q^{|A^{''}|}t)^{(-1)^{|A^{''}|}}\\
&=& \prod_{A^{''}} \prod _{S-A^{'} \subseteq A^{''}} \det \left( I-q^{|A^{''}|}t \alpha_{S-A^{'}, A^{'}}|B_{S-A^{'}, A^{'}} \right)^{(-1)^{|A^{''}|}} \\
&=& \prod_{A^{''}} \det \left( I-t (\alpha_{S-A^{'}, A^{'}})_{|A^{''}|}|K_{|A^{''}|}(S-A^{'}, A^{'}) \right)^{(-1)^{|A^{''}|}}
\end{eqnarray*}
where $S-A^{'} \subseteq A^{''}$ and $i=|A^{''}|$. Passing to homology, we have  
$$\prod_{A \subseteq S}L^*(x_0f_A, t)^{(-1)^n} = \prod_{A^{''}} \det \left( I-q^{|A^{''}|}t H_{|A^{''}|}(\alpha_{S-A^{'}, A^{'}})|H_{|A^{''}|}(K_*(S-A^{'}, A^{'})) \right)^{(-1)^{|A^{''}|}}$$
Since $f$ is $\bigtriangleup$-regular, so as $f_{S-A^{'}}$ for all $A^{'} \subseteq S$. By theorem 3.2, $H_{|A^{''}|}(S-A^{'}, \emptyset)=0$ for $|A^{''}|>0$. Lemma 4.6 implies that $H_{|A^{''}|}(S-A^{'}, A^{'})=0$ for $|A^{''}|>0$. For $|A^{''}|=0$, that means $A^{''}= \emptyset$ and $A^{'}=S$, we can conclude that 
\begin{eqnarray*}
\prod_{A \subseteq S}L^*(x_0f_A, t)^{(-1)^n} &=& \det \left( I-t H_0(\alpha_{\emptyset, S})|H_0(K_*(\emptyset, S)) \right)
\end{eqnarray*}
The zeta function of $\overline{U_f}$ can be expressed as
\begin{eqnarray*}
Z(\overline{U_f}, qt) = \prod_{\emptyset \neq \sigma \subseteq \bigtriangleup} Z(\gg_m^{\dim \sigma}, t) \prod_{\emptyset \neq \sigma \subseteq \bigtriangleup} L^*(x_0f_{\sigma}, t).
\end{eqnarray*}
Since 
$$\pp^n= \bigsqcup_{i=1}^{n+1}(\gg_m^{i-1})^{\binom{n+1}{i}}$$
where $i$ is the number of nonzero entries in the homogenuous coordinate in $\pp^n$, we have 
$$\sum_{\emptyset \neq \sigma \subseteq \bigtriangleup} N_k(\gg_m^{\dim \sigma})= \sum_{i=0}^n \binom{n+1}{i+1} (q^k-1)^i =N_k(\pp^n).$$
Then 
\begin{eqnarray*}
Z(\overline{U_f}, qt) &=& Z(\pp^n, t) \prod_{A \subsetneq S} L^*(x_0f_A, t)\\
&=& \frac{L^*(x_0f_S, t)^{-1}}{(1-t) \cdots (1-q^nt)} \prod_{A \subseteq S} L^*(x_0f_A, t).
\end{eqnarray*}
Substitute $A=S$ into lemma 4.3, we have 
\begin{eqnarray*}
L^*(x_0f_S, t)^{-1} &=& P_S(t)\\
&=& \det (I- t \alpha_S|K_*(S))\\
&=& 1-t
\end{eqnarray*}
where $\alpha_S$ is an identity endomorphism on $K_*(S)$. Hence  
\begin{eqnarray*}
Z(\overline{U_f}, qt) &=& \frac{\prod_{A \subseteq S} L^*(x_0f_A, t)}{(1-qt) \cdots (1-q^nt)} \\
&=& \frac{\det \left( I-t H_0(\alpha_{\emptyset, S})|H_0(K_*(\emptyset, S)) \right)^{(-1)^n}}{(1-qt) \cdots (1-q^nt)} . 
\end{eqnarray*}
$\hfill$ $\Box$\\

The mapping $H_0(\alpha_{\emptyset, S})$ is invertible as an endomorphism of $H_0(K_*(\emptyset, S))=B_{\emptyset, S}/$             $ \sum_{i=0}^n \widehat{D_i} B_{\emptyset, S- \{ i \}}$. It is easily seen that $H_0(\alpha_{\emptyset, S})$ is compatible with the action of Frobenius.\\

\noindent $\bf{Lemma \mbox{ 4.8.   }}$ Let $f \in \ff_q[x_1^{\pm 1}, \ldots , x_n^{\pm 1}]$ be $\bigtriangleup$-regular and suppose $\bigtriangleup (f)$ is an $n$-dimensional simplex and $p \nmid \det (a_{ij})$ where the matrix $(a_{i,j})$ is defined as in (5). For all $A \subseteq S$ and $A^{'} \subseteq S-A$, we have 
$$\det (I-t H_0(\alpha_{A, A^{'}})|H_0(K_*(A, A^{'})))^{(-1)^{|A|}}= \prod_{A \subseteq C \subseteq A \cup A^{'}}  \det (I-t H_0(\alpha_C)|H_0(K_*(C)))^{(-1)^{|C|}}$$ 
\\

\noindent $\bf{Proof.}$ Use induction on $|A^{'}|$ to the exact sequence (14). The proof is similar to the proof of lemma 4.4. $\hfill$ $\Box$\\ 

Hence, $\det \left( I-t H_0(\alpha_{\emptyset, S})|H_0(K_*(\emptyset, S)) \right)$ is a polynomial of degree 
\begin{eqnarray*}
\sum_{C \subseteq S} (-1)^{|C|} (n+1-|C|)! \mathrm{Vol}(x_0f_C).
\end{eqnarray*}

\section{$p$-adic estimates}
\   

In this section, we give a sharp lower bound for the Newton polygon of the polynomial $\det (I-t H_0(\alpha_{\emptyset, S})|H_0(K_*(\emptyset, S)))$.\\

Let $\bigtriangleup (f)$ be an $n$-dimensional simplex. For $A \subseteq S$, the graded ring $S_{\bigtriangleup}$ induces the graded subring  
$$S_{\bigtriangleup, A}= \left \{ \sum_{(r, m) \in M(x_0f)} A_{r, m} x_0^rx^m \in S_{\bigtriangleup} | A_{r, m}=0 \mbox{  if   }  (r, m) \notin \bigcap_{i \in A} H_i \right \}$$
with the grading $\deg (x_0^rx^m)=r$. Let $S_{\bigtriangleup, A}^r$ be a finite dimensional $\ff_q$-vector space with the basis $x_0^rx^m \in S_{\bigtriangleup, A}$. Define $F_{H_i}(x_0, x)= \sum_{j=0}^na_{ij}F_j(x_0, x) \in S_{\bigtriangleup}^1$ and $F_{H_i}^A(x_0, x)= \theta_A \circ F_{H_i}(x_0, x) \in S_{\bigtriangleup, A}^1$ where $(a_{ij})$ is defined as in (5).\\ 

\noindent $\bf{Lemma \mbox{ 5.1.   }}$ Let $f \in \ff_q[x_1^{\pm 1}, \ldots , x_n^{\pm 1}]$ be $\bigtriangleup$-regular and suppose $\bigtriangleup (f)$ is an $n$-dimensional simplex and $p \nmid \det (a_{ij})$ where the matrix $(a_{ij})$ is defined as in (5). For all $A \subseteq S$, the set $\{ F_{H_i}^A(x_0, x) \}_{i \notin A}$ taken in any order forms a regular sequence on $S_{\bigtriangleup, A}$, hence the Koszul complex $\overline{K}_*(S_{\bigtriangleup, A}, \{ F_{H_i}^A(x_0, x) \}_{i \notin A})$ is acyclic in positive dimension. Furthermore,  
$$\dim_{\ff_q} H_0( \overline{K}_*(S_{\bigtriangleup, A}, \{ F_{H_i}^A(x_0, x) \}_{i \notin A}))=(n+1-|A|)! \mathrm{Vol}(x_0f_A).$$

\noindent $\bf{Proof.}$ The hypotheses imply that $x_0f_A$ is $\bigtriangleup$-regular and $\bigtriangleup (x_0f_A)$ is also an $(n+1-|A|)$-dimensional simplex. The argument of Lemma 4.1 shows that one can simply repeat the proof of Theroem 3.1. $\hfill$ $\Box$\\ 

One defines $\theta_j:S_{\bigtriangleup, A} \rightarrow S_{\bigtriangleup, A \cup \{ j \}}$ as in (6). For $A^{'} \subseteq S-A$, this allows us to define, as in (10), subcomplexes $\overline{K}_*(A, A^{'})$ of the Koszul complex $\overline{K}_*(S_{\bigtriangleup, A}, \{ F_{H_i}^A(x_0, x) \}_{i \notin A})$. Proceeding as in the deduction of Lemma 4.8, we deduce:\\ 

\noindent $\bf{Lemma \mbox{ 5.2.   }}$ Let $f \in \ff_q[x_1^{\pm 1}, \ldots , x_n^{\pm 1}]$ be $\bigtriangleup$-regular and suppose $\bigtriangleup (f)$ is an $n$-dimensional simplex and $p \nmid \det (a_{ij})$ where the matrix $(a_{ij})$ is defined as in (5). For all $A \subseteq S$ and $A^{'} \subseteq S-A$, the complex $\overline{K}_*(A, A^{'})$ is acyclic in positive dimension. Furthermore,  
$$\dim_{\ff_q} H_0(\overline{K}_*(\emptyset, S))= \sum_{C \subseteq S} (-1)^{|C|} (n+1-|C|)! \mathrm{Vol}(x_0f_C).$$

Let $S_{\bigtriangleup, A, A^{'}}= \bigcap_{j \in A^{'}}(\ker \theta_j|S_{\bigtriangleup, A})$. By Lemma 5.2, $$H_0(\overline{K}_*(\emptyset, S))=S_{\bigtriangleup, \emptyset, S} / \sum_{i=0}^n F_{H_i}(x_0, x) S_{\bigtriangleup, \emptyset, S- \{ i \}}$$
is a finite dimensional graded ring. Let $h_{\bigtriangleup}(k)$ be the dimension over $\ff_q$ of graded piece of degree $k$ of $H_0(\overline{K}_*(\emptyset, S))$. It can be shown (as in \cite{Kou}, Lemma 2.9 and \cite{Dan}, Lemma 1.5) that $h_{\bigtriangleup}(k)=0$ for $i>n+1$. Define $v(\bigtriangleup)=\dim_{\ff_q} H_0(\overline{K}_*(\emptyset, S))$. Then $\sum_{k=0}^{n+1} h_{\bigtriangleup}(k)=v(\bigtriangleup)$.\\ 

Let $f$ be $\bigtriangleup$-regular over $\ff_q$, $\bigtriangleup (f)$ an $n$-dimensional simplex and $p \nmid \det (a_{ij})$ where the matrix $(a_{ij})$ is defined as in (5). Define the Hodge polygon $\mathrm{HP}(\bigtriangleup)$ in $\rr^2$ to be the polygon with the vertices $(0,0)$ and $(\sum_{k=0}^m h_{\bigtriangleup}(k),$ $\mbox{}$ $ \sum_{k=0}^m kh_{\bigtriangleup}(k))$ for $m=0, \ldots , n+1$. It can be shown (\cite{AS2}, Theorem 4.8) that $h(k)=h(n+1-k)$. Then  
$$\sum_{i=0}^{n+1}kh(k)= \frac{n+1}{2} v(\bigtriangleup).$$
\\

Write 
$$\det \left( I-t H_0(\alpha_{(\emptyset, S)})|H_0(K_*(\emptyset, S)) \right)= \sum_{m=0}^{v(\bigtriangleup)}A_mt^m,$$
\noindent where $A_0=1$ and $A_m \in \zz$. Define the Newton polygon $\mathrm{NP}(f)$ of $\sum_{m=0}^{v(\bigtriangleup)}A_mt^m$ to be the lower convex closure in $\rr^2$ of the points $(m, \mathrm{ord}_q(A_m))$ for $m=0, \ldots , v(\bigtriangleup)$.\\

\noindent $\bf{Theorem \mbox{   5.3.   }}$ Suppose $f \in \ff_q[x_1^{\pm 1}, \ldots , x_n^{\pm 1}]$ is $\bigtriangleup$-regular, $\bigtriangleup (f)$ is an $n$-dimensional simplex and $p \nmid \det (a_{ij})$ where the matrix $(a_{ij})$ is defined as in (5). Then the Newton polygon $\mathrm{NP}(f)$ lies above the Hodge polygon $\mathrm{HP}(\bigtriangleup)$ and their endpoints coincide.\\ 

\noindent $\bf{Proof.}$ This is proved by the method used in (\cite{AS1}, Theorem 3.10). The coincidence of the endpoints, which was not proved in \cite{AS1}, follows as in (\cite{Dwork2}, paragraph preceding Theorem 7.1). $\hfill$ $\Box$\\

\noindent $\bf{Remark \mbox{   5.4.    }}$ Let ${\cal{M}}_p(\bigtriangleup)$ be the moduli space of all $\bigtriangleup$-regular, Laurent polynomials over $\ff_p$ with $\bigtriangleup (f) =\bigtriangleup$ and ${\cal{H}}_p(\bigtriangleup)$ the moduli space of those $f \in {\cal{M}}_p(\bigtriangleup)$ such that its Newton polygon coincides with its Hodge polygon. It can be shown \cite{Wan} and \cite{Wan2} that there is an integer $D^*(\bigtriangleup)$ depending only on $\bigtriangleup$ such that if $p$ is a large prime and $p \equiv 1$ (mod $D^*(\bigtriangleup)$), then the Newton polygon of $\det \left( I-t H_0(\alpha_{(\emptyset, S)})|H_0(K_*(\emptyset, S)) \right)$ coincides generically with the Hodge polygon $HP(\bigtriangleup)$, i.e., for such $p$, ${\cal{H}}_p(\bigtriangleup)$ is a Zariski dense, open subset of ${\cal{M}}_p(\bigtriangleup)$.\\

\section{Calabi-Yau Hypersurfaces}
\   

Let $n \geq 2$ be a positive integer. We consider the universal family of Calabi-Yau complex hypersurfaces of degree $n+1$ in the projective space $\pp^n$. Its mirror family is a one-parameter family of toric hypersurfaces. To construct the mirror family, we consider the one-parameter subfamily $X_{\lambda}$ of complex projective hypersurfaces of degree $n+1$ in $\pp^n$ defined by  
$$f(x_1, \ldots , x_{n+1})= x_1^{n+1}+ \cdots + x_{n+1}^{n+1}+ \lambda x_1 \cdots x_{n+1}=0,$$
where $\lambda \in \cc$ is the parameter. The variety $X_{\lambda}$ is a Calabi-Yau manifold when $X_{\lambda}$ is smooth. Let $\mu_{n+1}$ denote the group of $(n+1)$-th roots of unity. Let 
$$G= \{ (\zeta_1, \ldots , \zeta_{n+1})| \zeta_i^{n+1}=1, \zeta_1 \ldots \zeta_{n+1}=1 \}/ \mu_{n+1} \cong ( \zz/(n+1) \zz)^{n-1},$$ 
where $\mu_{n+1}$ is embedded in $G$ via the diagonal embedding. The finite group $G$ acts on $X_{\lambda}$ by 
$$(\zeta_1, \ldots , \zeta_{n+1})(x_1, \ldots , x_{n+1})=(\zeta_1 x_1, \ldots , \zeta_{n+1}x_{n+1}).$$ 
The quotient $X_{\lambda}/G$ is a projective toric hypersurface $Y_{\lambda}$ in the toric variety $\pp_{\bigtriangleup}$, where $\bigtriangleup$ is the simplex in $\rr^n$ with vertices $\{ e_1 , \ldots , e_n, -(e_1+ \ldots +e_n) \}$ and the $e_i$'s are the standard coordinate vectors in $\rr^n$. Explicitly, the variety $Y_{\lambda}$ is the projective closure in $\pp_{\bigtriangleup}$ of the affine toric hypersurface in $\gg_m^n$ defined by
$$g(x_1, \ldots , x_n)=x_1+ \ldots +x_n+ \frac{1}{x_1 \ldots x_n}+ \lambda=0.$$
We are interested in the number of $\ff_q$-rational points on $Y_{\lambda}$.\\

The toric variety $\pp_{\bigtriangleup}$ has the following disjoint decomposition:
$$\pp_{\bigtriangleup}= \bigsqcup_{\sigma \in \bigtriangleup} \gg_m^{\dim \sigma} ,$$
where $\sigma$ runs over all non-empty faces of $\bigtriangleup$. Accordingly, the projective toric hypersurface $Y_{\lambda}$ has the corresponding disjoint decomposition 
$$Y_{\lambda}= \bigsqcup_{\sigma \in \bigtriangleup} Y_{\lambda, \sigma}, \mbox{                } Y_{\lambda, \sigma}=Y_{\lambda} \cap \gg_m^{\dim \sigma}.$$
For $\sigma = \bigtriangleup$, the subvariety $Y_{\lambda, \bigtriangleup}$ is simply the affine toric hypersurface defined by $g=0$ in $\gg_m^n$. For zero-dimensional $\sigma$, $Y_{\lambda, \sigma}$ is empty. For a face $\sigma$ with $1 \leq \dim \sigma \leq n-1$, one checks that $Y_{\lambda, \sigma}$ is isomorphic to the affine toric hypersurface in $\gg_m^{\dim \sigma}$ defined by 
$$1+x_1+ \cdots +x_{\dim \sigma} =0.$$
For such a $\sigma$, the inclusion-exclusion principle shows that
$$\# Y_{\lambda, \sigma}(\ff_q)=q^{\dim \sigma -1} - \binom{\dim \sigma}{1}q^{\dim \sigma -2}+ \cdots +(-1)^{\dim \sigma -1} \binom{\dim \sigma}{\dim \sigma -1} .$$
Thus,
$$\# Y_{\lambda, \sigma}(\ff_q)= \frac{1}{q} \left( (q-1)^{\dim \sigma} + (-1)^{\dim \sigma +1} \right) .$$
This formula holds even for zero dimensional $\sigma$ as both sides would then be zero.\\

Putting these calculations together, we deduce that
$$\# Y_{\lambda}(\ff_q)=\# Y_{\lambda, \bigtriangleup}(\ff_q) - \frac{(q-1)^n+(-1)^{n+1}}{q}+ \sum_{\sigma \in \bigtriangleup} \frac{1}{q} \left( (q-1)^{\dim \sigma} + (-1)^{\dim \sigma +1} \right) ,$$
where $\sigma$ runs over all non-empty faces of $\bigtriangleup$ including $\bigtriangleup$ itself. Since $\bigtriangleup$ is a simplex, one computes that   
$$\sum_{\sigma \in \bigtriangleup} \frac{1}{q} \left( (q-1)^{\dim \sigma} + (-1)^{\dim \sigma +1} \right) = \frac{q^{n+1}-1}{q-1}+(-1)= \frac{q(q^n-1)}{q-1}.$$
This implies that 
$$\# Y_{\lambda}(\ff_q)=\# Y_{\lambda, \bigtriangleup}(\ff_q) - \frac{(q-1)^n+(-1)^{n+1}}{q}+ \frac{q^n-1}{q-1}.$$
This equality holds for all $\lambda \in \ff_q$, including the case $\lambda =0$.\\

One of the basic problems in arithmetic mirror symmetry is to compare the number of rational points on a mirror pair $(X_{\lambda} , Y_{\lambda})$ of Calabi-Yau manifolds. Theorem 3.7 then gives us 
$$\# Y_{\lambda}(\ff_q)= \frac{q^n-1}{q-1} + \frac{(-1)^n \mathrm{Tr}(H_0(\alpha_{(\emptyset, S)})|H_0(K_*(\emptyset, S)))}{q}.$$
It can be shown that $\# X_{\lambda} (\ff_q) \equiv \# Y_{\lambda} (\ff_q)$ mod $q$ \cite{Wan1}.\\




\end{document}